\documentclass[reqno]{amsart} \usepackage{amscd}
\usepackage{epsf}
\newtheorem{theorem}{Theorem}[section]
\newtheorem{proposition}[theorem]{Proposition}
\newtheorem{lemma}[theorem]{Lemma}

\theoremstyle{remark} 
\numberwithin{equation}{section}
\newcommand{\field}[1]{\ensuremath{\mathbb{#1}}}
\newcommand{\CC}{\field{C}}
\newcommand{\HH}{\field{H}}
\newcommand{\RR}{\field{R}}
\newcommand{\ZZ}{\field{Z}}
\begin{document}

\title[Dimensional reductions]{Symplectic and Hyperk\"{a}hler structures in
a dimensional reduction of the Seiberg-Witten equations with a Higgs field}

\author{Rukmini Dey}

\begin{abstract}
In this paper we show that the dimensionally reduced Seiberg-Witten
equations lead to a Higgs field and  study the resulting moduli spaces.
The moduli space arising out of a subset of the equations, shown to be
non-empty for a compact Riemann surface of genus  $g \geq 1$,  gives rise to a
family  of moduli spaces carrying  a hyperk\"{a}hler structure.  For the full
set of equations the corresponding moduli space does not have the
aforementioned hyperk\"{a}hler structure but has a natural symplectic
structure. For the case of the torus,  $g=1$,  we show that  the full set of
equations   has a solution, different from the  ``vortex solutions''.
\end{abstract}

\maketitle

\section{Introduction}

Dimensional reductions of various gauge theories from four dimensions to
two dimensions have proved to be geometrically very rich. For example, the
moduli space of solutions to the dimensional reduction of the self-dual
Yang-Mills equations over a Riemann surface, ~\cite{H}, exhibits, among other
things, beautiful hyperk\"{a}hler structures.

It is important to study the analogous questions for  the Seiberg-Witten
equations. Though the main context  of Seiberg-Witten theory is in four
dimensions,  the $2$-dimensional reduction also seems worth exploring. A
reduction which gives the ``vortex'' equations have been studied extensively
e.g. Taubes ~\cite{T}, Bradlow, Garcia-Prada ~\cite{BG}, Olsen ~\cite{O},
Nergiz and Saclioglu ~\cite{NS}, ~\cite{NS2} and others. This reduction does
not have a Higgs field.

In the present paper we study a more general  dimensional reduction of  the
Seiberg-Witten equations  which gives three  equations. The novel feature is
the presence of   a Higgs field, which in our case is an imaginary valued
$1$-form. In the first few sections we consider the solutions to a subset of
the equations, namely, $[(2.1), (2.2)]$  for genus $g \geq 1$ compact Riemann
surfaces. The resulting moduli spaces which we denote by ${\mathcal M}$ and
$\Sigma_{\Psi}$ are hyperk\"{a}hler. This is partly due to the presence of
the Higgs field, $\Phi$.

In parallel we show that the full set of equations  $ [(2.1)-(2.3)]$ has
a solution with $\Phi \neq 0$ in the case of  genus $g=1$ compact Riemann
surface. We assume that solutions exist for genus $g >1$ as well. The
corresponding moduli space which we denote by ${\mathcal N}$ has a symplectic
and almost complex structure. The hyperk\"{a}hler structure, however, does
not descend to  ${\mathcal N}$. Setting $\Phi = 0$ results in ``vortex''
equations. We calculate the ``virtual'' dimensions of the moduli spaces.

\section{Dimensional Reductions of the Seiberg-Witten equations}

In this section we dimensionally reduce the Seiberg - Witten
equations on ${\RR }^4$ to  ${\RR }^2$ and define them over a
compact Riemann surface $M$.

\subsection{The Seiberg-Witten equations on ${\RR }^4$: }
This is a brief description of the Seiberg-Witten equations on
${\RR}^4$,~\cite{S},~\cite{Ak}, ~\cite{M} .

Identify ${\RR }^4 $  with the quaternions ${\HH}$ (coordinates $x = (x_1, x_2, x_3, x_4) $ identified with $\zeta=(\zeta_1, \zeta_2, \zeta_3, \zeta_4)$) and
let $\{e_i, i=1,2,3,4 \}$ be a  basis for ${\HH}$.
Fix the constant spin structure
$\Gamma : {\HH} = T_x {\HH} \rightarrow {\CC}^{4 \times 4}$, given by
$\Gamma (\zeta)  = \left[ \begin{array}{cc}
0              & \gamma (\zeta) \\
\gamma(\zeta) ^{*}  & 0
\end{array} \right], $
where $\gamma (\zeta) = \left[ \begin{array}{cc}
\zeta_1 + i \zeta_2   &-\zeta_3 - i \zeta_4 \\
\zeta_3 - i \zeta_4  & \zeta_1 - i \zeta_2
\end{array} \right]. $
 Thus $\gamma (e_1) = Id$, $\gamma (e_2) = I$, $\gamma (e_3) = J$,
$\gamma (e_4) = K$ where
$$I =  \left[ \begin{array}{cc}
i & 0 \\
0 & -i
\end{array} \right], J=  \left[ \begin{array}{cc}
0              & -1 \\
1             & 0
\end{array} \right], K =  \left[ \begin{array}{cc}
0              & -i \\
-i             &  0
\end{array} \right], $$
so that $I J = K$, $JK = I $, $KI = J$ and $I^2 = J^2 = K^2 = -Id$.

Recall that $Spin^{c}({\RR }^4) = (Spin ({\RR }^4) \times S^1) / {\ZZ}_2$.
Since $Spin ({\RR }^4)$ is a double cover of $SO(4)$, a $spin^c$ -
connection involves a connection $\omega$ on $T{\HH}$ and a connection
$ A = i\sum \limits_{j=1}^{4} A_j d x_j  \in \Omega^1 ({\HH}, i
{\RR})$ on
the characteristic line bundle $ {\HH} \times {\CC}$ which arises
from the $S^1$ factor (see ~\cite{S}, ~\cite{M}, ~\cite{Ak} for more
details). We set $\omega = 0$, which is equivalent to choosing the
covariant derivative  on the trivial tangent bundle to be $d$.  This is
legitimate since we are on ${\RR }^4$.   The  curvature $2$-form of the
connection $A$ is given by
$F(A) = d A    \in \Omega^2({\HH}, i{\RR})$. Consider the covariant
derivative acting on $\Psi \in C^{\infty}({\HH}, {\CC}^2)$ (the
positive spinor on ${\RR }^4$)  induced by the connection $A$ on
${\HH} \times {\CC}:$ $\nabla_j  \Psi = (\frac{\partial }{\partial
x_j} + i A_j) \Psi.  $ Then according to ~\cite{S}, the
Seiberg-Witten equations for $(A, \Psi)$ on ${\RR }^4$ are
equivalent to the equations:

$(SW1):$ $\nabla_1 \Psi = I \nabla_2
\Psi + J \nabla_3 \Psi+ K  \nabla_4 \Psi,$

$(SW2a):$  $F_{12} + F_{34} =  \frac{1}{2} \Psi^* I \Psi =
\frac{i}{2}(|\psi_1 |^2 - |\psi_2 |^2) \stackrel{\cdot}{=}\frac{1}{2}\eta_1,$

$(SW2b):$ $ F_{13} + F_{42} = \frac{1}{2} \Psi^* J \Psi = i (Im \psi_1 \bar{\psi}_2) \stackrel{\cdot}{=}\frac{1}{2} \eta_2, $

$(SW2c):$ $ F_{14} + F_{23} =
\frac{1}{2} \Psi^* K \Psi =  -i (Re \psi_1 \bar{\psi}_2)
\stackrel{\cdot}{=}\frac{1}{2} \eta_3 $

where $\Psi = \left[ \begin{array}{cc}
\psi_1 \\
\psi_2
\end{array} \right]$.

\subsection{Dimensional Reduction to ${\RR}^2$ }: Using the same method of
dimensional reduction as  in ~\cite{H}, we get the general form of the
reduced equations which contain the so-called Higgs field. Namely, impose the
condition that none of the $A_i$'s and  $\Psi$ in $(SW1)$ and $(SW2)$ depend
on $x_3$ and $x_4$, i.e. $ A_i = A_i(x_1, x_2) $, $\Psi = \Psi(x_1, x_2)$ and
 set $\phi_1 = -i A_3$ and $\phi_2 = -i A_4$. The $(SW2)$
 equations reduce to the following system on ${\RR }^2$, $ F_{12}
=\frac{1}{2} \eta_1, $ and two other equations, which, after
 introducing complex coordinates $z=x_1 + i x_2 $, can be rewritten
 as: $ \frac{\partial (\phi_1 + i \phi_2 )}{ \partial \bar{z} } =
 -\frac{1}{2}(\eta_2 + i\eta_3) = -\psi_1 \bar{\psi}_2 , $ where
 $\frac{\partial}{\partial \bar{z}} =
\frac{1}{2}(\frac{\partial}{\partial x_1} + i
\frac{\partial}{\partial x_2 }). $ Setting $ \phi_1 + i \phi_2 =
\phi $ and $\omega = i dz \wedge d \bar{z}$ we  rewrite the reduction of
$(SW2)$ as the following two equations,

$(1)$ $F(A) = \frac{i}{2} ( |\psi_1|^2 - |\psi_2|^2 ) \omega, $

$(2)$ $ 2\bar{\partial} \Phi = -i (\psi_1 \bar{ \psi_2} ) \omega $

where
$ \Phi = \phi dz - \bar{\phi} d \bar{z} \in \Omega^1 ({\RR }^2,
{i\RR }) \rm{\; and \;}  \psi_1 , \psi_2  \in C^{\infty} ({\RR}^2, {\CC}) $
are spinors on ${\RR }^2$. Next  consider the Dirac
equation   $(SW1)$:

$ \nabla_1 \psi - I \nabla_2 \psi - J \nabla_3 \psi - K \nabla_4 \psi = 0 $
which is rewritten as

$\left[ \begin{array}{cc}
\frac{\partial}{\partial x_1 } + iA_1 - i \frac{\partial}{\partial x_2} + A_2
 & \frac{\partial}{\partial x_3} + iA_3 + i \frac{\partial}{\partial x_4} -
A_4     \\
-\frac{\partial}{\partial x_3 } - iA_3 + i \frac{\partial}{\partial x_4} -
A_4  & \frac{\partial}{\partial x_1} + iA_1 +  i \frac{\partial}{\partial
x_2} - A_2 \end{array}
\right] \left[ \begin{array}{cc}
\psi_1 \\
\psi_2
\end{array} \right]   = 0.$

Introducing $A = \frac{i}{2}(A_1 - i A_2) dz$ and
$-\bar{A} = \frac{i}{2} (A_1 + i A_2) d \bar{z}$ where the total connection
$A - \bar{A} =i (A_1 dx + A_2 dy)$ we can finally write it as

$(3)$ $ \left[ \begin{array}{cc} -\frac{1}{2}\bar{\phi}d \bar{z} &
(\bar{\partial} -\bar{A} )  \\
(\partial + A)  & -\frac{1}{2}\phi dz
\end{array} \right] \left[ \begin{array}{cc}
 \psi_1 \\
\psi_2
\end{array} \right] = 0  $

We call equations $(1) - (3)$ as the dimensionally reduced
Seiberg-Witten equations over ${\CC}$.

\subsection{The Dimensionally Reduced Equations on a Riemann surface}

Let $M$ be a compact Riemann surface of genus $g \geq 1$ with a
conformal metric $ds^2 = h^2 dz \otimes d \bar{z}$ and let $\omega
= i h^2 d z \wedge d \bar{z}$ be a real form proportional to the
induced K\"{a}hler form. Let  $L$ be a  line bundle with a
Hermitian metric $H$.  Let $\psi_1,\psi_2$ be sections of the line
bundle $L$ i.e., $\psi_1, \psi_2 \in \Gamma(M,L)$. Then we have an
inner product $ <\psi_1 , \psi_2 >_H $ and  norm $|\psi |_H \in
C^{\infty} $ of the sections of $L$. Let $A -\bar{A}$ be a unitary
connection on $L$ and $\Phi = \phi dz -\bar{\phi} d \bar{z} \in
\Omega^1 (M , i{\RR })$. We will  assume that $\Psi = \left[
\begin{array}{c}
\psi_1\\
\psi_2
\end{array} \right]$ is not identically zero. We can
rewrite the equations $(1) - (3)$ in an invariant form on $M$ as
follows:

 $$F(A) = i \frac{( |\psi_1 |_H ^2 - |\psi_2 |_H ^2 )}{2}
\omega, \leqno{(2.1)}$$

$$2\bar{\partial} \Phi  =  - i <\psi_1, \psi_2>_H
\omega, \leqno{(2.2)}$$

$$\left[ \begin{array}{cc} -\frac{1}{2} \bar{\phi}d \bar{z} &
\bar{\partial} -  \bar{A}  \\
\partial + A & -\frac{1}{2} \phi dz
\end{array} \right] \left[ \begin{array}{cc}
\psi_1 \\
\psi_2
\end{array} \right] = 0.\leqno{(2.3)}$$

{\bf Note:}  Equation $(2.3)$ has two equations. One equation $(2.3a)$ comes
with $-\bar{A}$ and the other one $(2.3 b)$  with $A$. This is unitarity of
the connection $A - \bar{A}$ on the line bundle, ~\cite{GH}. Setting
$\Phi = 0$ in equation $(2.2)$ and $(2.3)$ we obtain the usual vortex
equations (where either $\psi_1$ or $\psi_2$ is zero).

Let ${\mathcal C} = {\mathcal A} \times \Gamma (M, L \oplus L)
\times {\mathcal H}$ , where ${\mathcal A}$ is the space of
connections on a line bundle $L$, $\Gamma (M, L \oplus L) $ the
space of sections of the  bundle $L \oplus L$ and ${\mathcal H}$
be $\Omega^{1}(M, i {\RR})$, the space of Higgs fields. The gauge
group ${\mathcal G} = Maps (M, U(1))$ acts on ${\mathcal B}$ as
$(A, \Psi, \Phi) \rightarrow (A + u^{-1} du, u^{-1} \Psi, \Phi)$
and leaves the space of solutions to $(2.1) - (2.3) $ invariant.
There are no fixed points of this action. Because  a fixed point
would mean that there is a connection $A_0$ such that $A_0 +
u^{-1}du = A_0$ for all $u$ in the gauge group. This is not
possible. We assume throughout that $\Psi$ is not identically
zero.

Taking the quotient by the gauge group of the solutions to $(2.1) - (2.3)$ we
obtain a moduli space which we denote by ${\mathcal N}$. Let us denote the
moduli space of solutions to $(2.1)-(2.2)$  as ${\mathcal M}$ where we let
the equivalence class of  $\Psi$ vary. We define a new moduli space
$ \Sigma_{[\Psi]}$, by fixing an equivalence class of $\Psi $ as  follows.
Choose an appropriate $\Psi$ such that $< \psi_1, \psi_2> \omega $ is
$\bar{\partial}$-exact and let ${\mathcal W} \stackrel{\cdot}{=} {\mathcal A} \times \{ G \cdot \Psi \} \times {\mathcal H}$ $\subset$ ${\mathcal E}$ where
$\{ G \cdot \Psi \}$ is the orbit of $\Psi$ due to action of the gauge group.
Let ${\mathcal S}_1 \stackrel{\cdot}{=} {\mathcal W} \cap \tilde{S}, $ where
$\tilde{S}$ is the solution space to equations $(2.1)$ and $(2.2)$ on
${\mathcal C}$. Define $\Sigma_{[\Psi]} \stackrel{\cdot}{=} {\mathcal S}_1 / G $. Any point $p $ $\in$ $ \Sigma_{[\Psi]}$ is given by $p = ([(A, \Psi, \Phi])$ where $\Psi$ is now fixed, $[\cdot, \cdot, \cdot]$ denotes the gauge
equivalence class and $(A, \Psi, \Phi)$ satisfy equations $(2.1)$ and $(2.2)$.
 ${\mathcal S}_1 $ essentially consists of
$(A, \Phi)\in {\mathcal A} \times {\mathcal H}$ such that $d A = \frac{i}{2} f_1 \omega$ and $2\bar{\partial} \Phi = f_2 \omega,$ where $f_1 = | \psi_1 |^2_H - |\psi_2|^2_H$ and $f_2 =-i < \psi_1, \psi_2 >_H$  $\in C^{\infty}(M)$. We will see that if $(A, \Phi) \in \Sigma_{[\Psi]}$, then if one changes
$A \rightarrow A^{\prime} = A + \alpha$ such that
$d \alpha =0$, $\alpha$ unique upto exact forms, and $\Phi^{0,1} \rightarrow \Phi^{\prime 0,1 }= \Phi^{0,1} + \eta^{0,1}$ such that $\bar{\partial} \eta^{0,1} = 0$ then, the $(A^{\prime}, \Phi^{\prime}) \in \Sigma_{[\Psi]}$. Thus $\Sigma_{[\Psi]}$ is an affine space.

\begin{proposition}
If $L$ is a trivial line bundle on a compact Riemann surface
of genus $g=1$ then $(2.1)-(2.2)$ have a solution ( with $\Psi \neq 0$,
$\Phi \neq 0$).  If $L=K^{-1}$ on a compact Riemann surface of genus $>1$ then
$(2.1)-(2.2)$ has a solution (with $\Psi \neq 0$, $\Phi \neq 0$).
Thus ${\mathcal M}$, $\Sigma_{\Psi}$ is non-empty.
\end{proposition}

\begin{proof}
For genus $g=1$ let us take a metric of the form $ds^2 = dz \otimes d \bar{z}$.
Let $\Psi_1 = 1$ and $\Psi_2 = e^{i\beta} = e^{i(kz + \bar{k}\bar{z})}$ be an
eigenfunction of the Laplacian with eigenvalue $-1$, i.e. $|k| = 1$.
Let $A$ be any flat connection and
$\phi dz  = \frac{\partial (e^{-i \beta})}{2}$. They satisfy
$(2.1)$ and $(2.2)$. The reason for this kind of solution will be clear in
proposition (\ref{solution}).

For $g>1$, $ L =K^{-1}$ has a metric same as the metric on the surface
$ds^2 = h^2 dz \otimes d \bar{z}$, so that we can write $h$ instead of $H$.
We take $A - \bar{A}$ to be the usual connection induced by the metric i.e.
$A = \frac{\partial log h}{\partial z}$ and $\bar{A} =
 \frac{\partial log h}{\partial \bar{z}} $ ~\cite{GH}.
Then $F(A) = -(\frac{\partial ^2 log h}{\partial z \partial \bar{z}})
d z \wedge d \bar{z} = -i K \omega$ where
$K = \frac{- \Delta_h log h}{h^2}$ is the Gaussian curvature of the metric
$h$ and $\omega = \frac{i}{2} h^2 dz \wedge d \bar{z}$ is the K\"{a}hler
form corresponding to the metric. Thus the equation $(2.1)$ just reduces to
$K= \frac{|\psi_2|_h^2 - | \psi_1|^2_h }{2}$. The right hand side is a smooth
function on the Riemann surface. We know that as long as $|\psi_2|_h^2 <
|\psi_1|_h^2$, there is a solution when the genus is $>1$. In fact 
 there exists a metric, in every conformal class, such
that any arbitrary negative definite function can be admitted as a
Gaussian curvature of a Riemann surface of genus $>1$ ~\cite{B},
~\cite{De}, ~\cite{KW}.

Let $\Phi = \partial w$. Also let $\psi_1 = v \psi_2$, where $v$
is a function. This is possible since $\psi_1, \psi_2$ are both
sections of the same line bundle. Then equation $(2.2)$ becomes
$\Delta_h w = \frac{1}{h^2} \frac{\partial^2 w}{ \partial z \partial \bar{z}}
= \frac{-\bar{v}}{2} |\psi_2|_h^2 = \tau$ where $\Delta_h$ is the
Laplacian induced by the metric $h$. We choose $\psi_2$
arbitrarily. Now we choose $v$ to be such that
$\int_M \frac{-\bar{v}}{2} |\psi_2|_h^2 h^2 dz d \bar{z}= 0$ and
$|\psi_2|_h^2 < |\psi_1|_h^2$ hold. By Hodge theory ~\cite{GH} there exists a
Green's operator $G$ for the Laplacian such that $w= G\tau$ is a solution to
equation $(2.2)$.
\end{proof}

\begin{proposition}
Let $L$ be a trivial line bundle  on a compact Riemann surface of
genus $g = 1$ then $(2.1)-(2.3)$ have a solution with
$\Psi \neq 0, \Phi \neq 0$. Thus  ${\mathcal N}$ is non-empty. \label{solution}
\end{proposition}

\begin{proof}
Let us solve for the case of the case of the torus, $g=1$. 
Let our torus  be thought of as $0 \leq x \leq 2 \pi$ and
$0 \leq y \leq  2 \pi$ with the endpoints identified. We take the metric on the torus to be $ds^2 = dz \otimes d \bar{z}$, i.e. $h=1$.
The equations are then as follows

$(2.1)$ $ F(A) = - \frac{|\Psi_1|^2 - |\Psi_2|^2}{2} dz \wedge d\bar{z} =0 $

$(2.2)$ $ \bar{\partial} \Phi = \frac{-1}{2} \Psi_1 \bar{\Psi}_2  d \bar{z} \wedge d z $

$(2.3a)$ $\frac{\bar{\partial} \Psi_2}{\Psi_2} -  \bar{A} -\frac{1}{2} (\bar{\phi} d \bar{z}) \frac{\Psi_1}{\Psi_2} = 0$

$(2.3b)$ $\frac{\partial \Psi_1}{\Psi_1} +  A - \frac{1}{2} \phi dz \frac{\Psi_2}{\Psi_1} = 0$.

where $\Phi = \phi dz - \bar{\phi} d \bar{z}$

Since we took the line bundle to be trivial, one solution would be to take 
$\Psi_1= c_1$, $\Psi_2 = c_1 e^{i c_2(z + \bar{z})}$, $\phi dz = -i c_2 e^{-i c_2 (z + \bar{z})} dz $, $A= -\frac{i c_2}{2} dz$ where $c_1$ is a  complex
constant  and  $c_2$ is a real constant  satisfying  $|c_1| = \sqrt{2} c_2$. 

Near this solution there is a $4$-dimensional moduli space, see proposition
(~\ref{dimvor}).
\end{proof}

\begin{proposition}
Let us consider the moduli spaces $\Sigma_{\Psi}$, ${\mathcal M}$,
${\mathcal N}$, respectively.
 Suppose $(A, \Psi, \Phi)$ is a point on the moduli space such
that $\Psi$ is not identically  $0$. The (virtual) dimension of
$\Sigma_{\Psi}$ is $4g$, and  ${\mathcal M}$ is infinite
dimensional. The (virtual) dimension  of ${\mathcal N}$ is $2 g +
2$. If either $\Phi =0$ and $\psi_1 $ or $\psi_2$ is zero (the
vortex case) then the dimension of ${\mathcal N}$ is $c_1(L) + g +
1$ or $c_1(\bar{L}) +g + 1$, respectively.  \label{dimvor}
\end{proposition}

\begin{proof}
To calculate the dimension of $\Sigma_{[\Psi]}$, we linearize equations
$(2.1)$ and $(2.2)$ with equivalence class of $\Psi$ fixed to obtain:
\begin{eqnarray*}
 (I) & & d \alpha = 0 \\
(II) & & \bar{\partial} \eta^{1,0} = 0,
\end{eqnarray*}
where $(\alpha, \beta, \eta) \in T_{p} {\mathcal W}$ and $ \eta = \eta^{1,0} +
\eta^{0,1}$. Taking into account the gauge group action,  we get
dim[$\{ \alpha \in \Omega^{1} (M, i{\RR }) | d \alpha = 0 \} / \{ \alpha = df \} ] = 2g$. Also, dim[$\{ \eta \in \Omega^{1,0} (M, {\CC})| \bar{\partial} \eta = 0 \}]= 2g.$ Thus  the dim[$T_p \Sigma_{[\Psi]}]= 4g$.

${\mathcal M}$ is infinite dimensional since $\Psi$ is not fixed.

To calculate the dimension of ${\mathcal N}$  let ${\mathcal S}$ be the
solution space to $(2.1)-(2.3)$.
Consider the tangent space $T_p {\mathcal S}$ at a point
  $p= (A, \Psi, \Phi) \in {\mathcal S},$ which is defined by the
 linearization of equations $(2.1)-(2.3)$. Let $X = (\alpha, \beta,
\gamma )$ $\in T_p {\mathcal S} $, where $\alpha \in \Omega^{1}
(M, i {\RR })$ and $ \beta =\left[ \begin{array}{cc}
\beta_1 \\
\beta_2
\end{array} \right]  \in \Gamma (M, L \oplus L), $ and $\gamma \in {\mathcal H}$. The
linearizations of the equations are as follows

$$d\alpha = \frac{i}{2}( \beta_1 \bar{\psi_1} + \psi_1 \bar{\beta_1} - \beta_2 \bar{\psi_2} - \psi_2 \bar{\beta_2})  \omega, \leqno{(2.1)^{\prime}}$$

$$\bar{\partial} \eta^{1,0} = -i \frac{1}{2}
(\psi_1 \bar{\beta_2} + \beta_1 \bar{\psi_2}) \omega, \leqno{(2.2)^{\prime}}$$

$$ \left[ \begin{array}{cc}
-\frac{1}{2} \bar{\phi} \bar{dz}   & (\bar{\partial} - \bar{A}) \\
\partial + A  & -\frac{1}{2} \phi dz
\end{array} \right] \left[
\begin{array}{c}
\beta_1 \\
\beta_2
\end{array} \right]    +  \left[ \begin{array}{cc}
-\frac{1}{2} \gamma^{0,1}  & -\bar{\alpha}   \\
\alpha  & -\frac{1}{2} \gamma^{1,0}
\end{array} \right] \left[
\begin{array}{cc}
\psi_1 \\
\psi_2
\end{array} \right] = 0. \leqno{(2.3)^{\prime}} $$

Taking into account the  quotient by the gauge group ${\mathcal  G}$, we
arrive at the following  sequence ${\mathcal C}$
$$0 \rightarrow \Omega^0(M, i{\RR }) \stackrel{d_1}{\rightarrow}
\Omega^1 (M, i{\RR }) \oplus \Gamma (M, {\mathcal L}) \oplus {\mathcal H}
\stackrel{d_2}{\rightarrow} \Omega^2 (M, i{\RR }) \oplus \Omega^{2}(M, {\CC})
\oplus    V \rightarrow 0 ,$$
where  ${\mathcal L} = L \oplus L$, $V =  (L \otimes \Omega^{0,1}(M) ) \oplus  (L \otimes \Omega^{1,0}(M)) $,

$d_1 f = (df, -f \Psi, 0)$, $ d_2 (\alpha, \left[ \begin{array}{c}
\beta_1 \\
\beta_2
\end{array} \right] , \gamma)   \stackrel{\cdot}{=} (A, B, C ), $

$A= d\alpha - \frac{i}{2} [(\psi_1 \bar{\beta}_1 +
\beta_1 \bar{\psi}_1) - (\psi_2 \bar{\beta}_2 +
\beta_2 \bar{\psi}_2)] \omega \in \Omega^{2}(M, i {\RR})$

$B= \bar{\partial} \gamma^{1,0} + \frac{i}{2} (\psi_1
\bar{\beta}_2 + \beta_1 \bar{\psi}_2) \omega \in \Omega^{2}(M,{\CC})$

$C= \left[ \begin{array}{cc}
-\frac{1}{2} \bar{\phi} \bar{dz}   & (\bar{\partial} - \bar{A})   \\
(\partial + A)  & -\frac{1}{2} \phi dz
\end{array} \right] \left[\begin{array}{c}
 \beta_1 \\
\beta_2
\end{array} \right] +  \left[
\begin{array}{cc}
-\frac{1}{2} \gamma^{0,1} & -\bar{\alpha} \\
\alpha  & -\frac{1}{2} \gamma^{1,0} \end{array}
\right] \left[
\begin{array}{c}
\psi_1 \\
\psi_2
\end{array} \right] \in V.$

It is easy to check that $d_2 d_1 = 0$, so that this is a complex.
Clearly, $H^0({\mathcal C})$ $=$ $0$ ,  because if  $f \in $
 $ker(d_1),$ then $df =0$ and $f \Psi =0$, which  implies  $f = 0$ since we
are in the neighbourhood of a point where $\Psi \neq 0$.

The Zariski dimension of the moduli space is $dim H^1({\mathcal C})$
while the virtual dimension is
dim $H^1({\mathcal C}) -$ dim $H^2({\mathcal C})$,
and coincides with the Zariski dimension whenever $dim H^2 ({\mathcal C})$ is
zero (namely the smooth points of the solution space ~\cite{M}, page 66).
The virtual dimension is $ = dim H^1({\mathcal C}) - dim H^2({\mathcal C}) =$
index of ${\mathcal C} $.

To calculate the index of ${\mathcal C}$ , we consider the family
of complexes $({\mathcal C}^t, d^t)$, $0 \leq t \leq 1,$ where

$d_{1}^t = (df, -t f \Psi, 0)$,
$d_2 (\alpha, \left[ \begin{array}{cc}
\beta_1 \\
\beta_2
\end{array} \right] , \gamma)   \stackrel{\cdot}{=} (A_t, B_t, C_t)$,

$A_t =d\alpha - \frac{it}{2} [(\psi_1 \bar{\beta}_1
+ \beta_1 \bar{\psi}_1) - (\psi_2 \bar{\beta}_2 + \beta_2 \bar{\psi}_2)]
\omega $,

$B_t = \bar{\partial} \gamma^{1,0} + \frac{it}{2} (\psi_1 \bar{\beta}_2  + \beta_1 \bar{\psi}_2) \omega $,

$C_t= \left[ \begin{array}{cc}
-\frac{t}{2} \bar{\phi} \bar{dz}   & (\bar{\partial} - \bar{A})   \\
\partial + A  & -\frac{t}{2} \phi dz
\end{array} \right] \left[ \begin{array}{cc}
\beta_1 \\
\beta_2
\end{array} \right]  + t \left[ \begin{array}{cc}
                                      -\frac{1}{2} \gamma^{0,1} &  -\bar{\alpha} \\
\alpha & - \frac{1}{2} \gamma^{1,0}
                                        \end{array} \right] \left[
                                       \begin{array}{c}
                                       \psi_1 \\
                                       \psi_2
                                       \end{array} \right] )$.

Clearly, $ind ({\mathcal C}^t)$ does not depend on $t$. The complex
${\mathcal C}^0$ (for $t=0$ ) is
\begin{eqnarray*}
 0 \rightarrow \Omega^0(M, i{\RR })\stackrel{d^{\prime}_1}{\rightarrow}
\Omega^1 (M, i{\RR }) \oplus \Gamma (M, {\mathcal L}) \oplus {\mathcal H}
\stackrel{d^{\prime}_2}{\rightarrow} \Omega^2 (M, i{\RR }) \\
\oplus \Omega^{2}(M, {\CC}) \oplus V \rightarrow 0
\end{eqnarray*}
\begin{eqnarray*}
\rm{where \;} & & d_1^{\prime} f = (df, 0, 0), d_2^{\prime} (\alpha, \beta, \gamma) \\
&=& (d\alpha, \bar{\partial} \gamma^{1,0},{\mathcal D}_A \beta).
\end{eqnarray*}
 Here  ${\mathcal D}_A =   \left[ \begin{array}{cc}
0  & \bar{\partial} - \bar{A}    \\
\partial + A & 0
\end{array} \right]. $

${\mathcal C}^0$ decomposes into a direct sum of three complexes

$(a)$ $ 0 \rightarrow \Omega^0(X, i{\RR }) \stackrel{ d}{\rightarrow} \Omega^1 (M, i{\RR }) \stackrel{d}{\rightarrow}  \Omega^2 (M, i{\RR }) \rightarrow 0,$

$(b)$   $ 0 \rightarrow \Omega^{1}(M, i {\RR}) \stackrel{\bar{\partial}}{\rightarrow} \Omega^{1,1} (M,i{\RR}) \rightarrow 0 $

$(c)$  $ 0 \rightarrow \Gamma (M,S) \stackrel{{\mathcal
D}_A}{\rightarrow} \Gamma (M,S^{\prime}) \rightarrow 0,$ where $S=
L \oplus L$, $S^{\prime} = (L \otimes K) \oplus (L \otimes
\bar{K})$.

dim $H^1$(complex (a))$=2g$, dim $H^1$(complex (b))$=2g$.

The complex $(c)$ breaks into two complexes as follows

  $ (c1) $ $ 0 \rightarrow \Gamma (M,L)
\stackrel{\partial + A}{\rightarrow} \Gamma (M, L \otimes K)
\rightarrow 0. $

$ (c2) $ $ 0 \rightarrow \Gamma (M, L)
\stackrel{\bar{\partial} - \bar{A}}{\rightarrow}\Gamma (M, L \otimes \bar{K}) \rightarrow 0. $

$(c1)$ comes from the equation $(\partial + A) \psi_1  = 0$.
Taking complex conjugate one gets $(\bar{\partial} + \bar{A}) \bar{\psi}_1  = 0$ which is the holomorphicity of a section of $\bar{L}$. Thus the first complex in $(c)$ can be rewritten as

$(c1)$ $0 \rightarrow \Gamma (M, \bar{L})
\stackrel{\bar{\partial} + \bar{A}}{\rightarrow} \Gamma (M,\bar{L}
\otimes \bar{K}) \rightarrow 0. $

By Riemann Roch, the index of $(c1)$ is $(c_1 ( \bar{L}) - g + 1)$ and
that of $(c2)$ is $(c_1 ( L ) - g + 1) $ and  thus the sum
is $2g + 2g + c_1(\bar{L}) -g + 1 + c_1(L) - g + 1$ or $2g+ 2$.
If $\psi_1 =0$ then the dimension is $2g + c_1(L) - g + 1 = c_1(L) + g + 1.$
If $\psi_2 =0$ then the dimension is $c_1(\bar{L}) + g + 1$.
\end{proof}

{\bf Note:} For simplicity of the exposition (to avoid writing the norms
explicitly in terms of the metric $H$) we choose $\bar{L} = L^{-1}$ ,
so that $\| \psi_1 \|^2 ,  \| \psi_2 \|^2 , <\psi_1, \psi_2>  \in
C^{\infty} (M, {\CC})$ are well defined. Thus throughout the rest
of the paper, we shall drop the subscripts $H$ in equation $(2.1)-(2.3)$.
We must mention that this assumption is not essential, but it is just to make
the exposition simpler.

\section{Symplectic and almost complex structures}
In the next theorem we discuss  symplectic and   complex
structures  on ${\mathcal N}$. For  similar work
on the vortex moduli space, see ~\cite{Br}, ~\cite{GP}.

Let ${\mathcal C} = {\mathcal A} \times \Gamma (M, L \oplus L)
\times {\mathcal H}$ be the space on which equations $(2.1) -
(2.3)$ are imposed. Let $p= (A, \Psi, \Phi) \in {\mathcal C}$, $X
= ( \alpha_1, \beta, \gamma_1)$, $Y= (\alpha_2, \eta, \gamma_2)$
$\in T_p {\mathcal C}$. On ${\mathcal C}$ one can define a metric
\begin{eqnarray*}
 g( X, Y) = \int_M *\alpha_1 \wedge \alpha_2  +  \int_M Re < \beta, \eta> \omega + \int_M * \gamma_1 \wedge \gamma_2
\end{eqnarray*}
and an almost complex structure  ${\mathcal I} = \left[
\begin{array}{cccc}
* & 0  & 0 & 0\\
0 & i & 0  & 0\\
0 & 0  & -i & 0\\
0 & 0 & 0 & *
\end{array} \right] : T_p {\mathcal C} \rightarrow T_p {\mathcal C}$
where   $*: \Omega^{1} \rightarrow \Omega^{1}$ is  the Hodge star
operator on $M$ (which takes type $dx$ forms to type $dy$ and $dy$
to $-dx$ , i.e $*(\eta dz) = -i \eta  dz, *(\eta d \bar{z})= i
\eta d \bar{z}$ ). We define
\begin{eqnarray*}
\Omega(X, Y) = -\int_{M} \alpha_1 \wedge \alpha_2 + \int_{M} Re <
I \beta , \eta> \omega - \int_M \gamma_1 \wedge \gamma_2
\end{eqnarray*}
where $ I = \left[\begin{array}{cc}
i & 0 \\
0  & -i
\end{array} \right]$
such that $ g ({\mathcal I} X, Y) = \Omega ( X, Y).$
Moreover, we have the following:

\begin{proposition} The metrics $g$,  the symplectic form $\Omega$,
and the almost complex structure  ${\mathcal I}$ are invariant
under the gauge group action on ${\mathcal C}$. \label{inv}
\end{proposition}

\begin{proof}
 Let $p = (A, \Psi, \Phi) \in {\mathcal C}$ and $u \in G, $ where
$u \cdot p = (A + u^{-1} du, u^{-1} \Psi, \Phi)$.

Then $u_* : T_p {\mathcal C} \rightarrow T_ {u \cdot p} {\mathcal C} $ is given
by the mapping $(Id, u^{-1}, Id)$ and it is now easy to check that $g$ and
$\Omega$ are invariant and ${\mathcal I}$ commutes with $u_*$.
\end{proof}

\begin{proposition}
\label{props} The equation $(2.1)$  can be realised as a moment
map $\mu = 0$ with respect to the action of the gauge group and
the symplectic form $\Omega$. \label{moment1}
\end{proposition}

\begin{proof}
Let $\zeta \in \Omega(M, i{\RR })  $ be the Lie algebra of the
gauge group (the gauge group element being $u = e^{ \zeta}$ );
It generates a vector field $X_{\zeta}$ on ${\mathcal C}$ as follows :
$$X_{\zeta} (A, \Psi, \Phi) = (d \zeta, -\zeta \Psi, 0) \in T_p
{\mathcal C},p = (A, \Psi, \Phi) \in {\mathcal C}.$$

We show next that $X_{\zeta}$ is Hamiltonian. Namely, define
$H_{\zeta} : {\mathcal C} \rightarrow {\CC} $ as follows: $$
H_{\zeta} (p) = \int_{M} \zeta \cdot (  F_{A} - i\frac{( |\psi_1|
^2 - |\psi_2| ^2 )}{2} \omega). $$  Then for $X = (\alpha,
\beta, \gamma) \in T_p {\mathcal C}$.
 \begin{eqnarray*}
 dH_{\zeta} ( X ) & = & \int_M \zeta d \alpha  -i \int_M \zeta  Re ( \psi_1 \bar{\beta_1} - \psi_2 \bar{ \beta_2} )  \omega    \\
 &= &\int_M (-d \zeta) \wedge \alpha  -  \int_M Re < I \zeta
\left[ \begin{array}{cc}
\psi_1 \\
\psi_2
\end{array} \right] , \left[ \begin{array}{cc}
\beta_1 \\
\beta_2
\end{array} \right] > \omega \\
 & = & \Omega ( X_{\zeta},  X ),
 \end{eqnarray*}
where we use that $\bar{\zeta} = - \zeta$.

 Thus we can define the moment map $ \mu : {\mathcal C} \rightarrow
 \Omega^2 ( M, i{\RR} )= {\mathcal G}^* $ ( the dual of the Lie
 algebra of the gauge group)  to be $$ \mu ( A, \Psi)
 \stackrel{\cdot}{=} (F(A) - i\frac{( | \psi_1 | ^2 -
 |\psi_2|^2)}{2}  \omega). $$ Thus equation $(2.1))$ is $\mu = 0$.
 \end{proof}

 \begin{lemma}
 Let $S$ be the solution spaces to equation $(2.1)- (2.3)$,  $X \in
 $  $ T_p {\mathcal S}$. Then ${\mathcal I}X $ $\in T_p {\mathcal S}$
 if and only if $X$ is orthogonal to the gauge orbit $ O_p = G
 \cdot p$. \label{ortho}
 \end{lemma}

 \begin{proof}
 Let $X_{\zeta} \in T_p O_p,$ where $\zeta $ $\in $ $\Omega^{0} (M,
+i {\RR })$,
$g( X , X_{\zeta} ) = -\Omega ({\mathcal I} X, X_{\zeta} ) =
- \int_M \zeta \cdot d \mu ( {\mathcal I} X ),$ and therefore
${\mathcal I} X$ satisfies the linearization of equation $(2.1)$ iff
$ d \mu ({\mathcal I} X)  = 0$, i.e.,  iff $g (X, X_{\zeta}) = 0$ for all
$\zeta$. Second, it is easy to check that ${\mathcal I} X$ satisfies the
 linearization of equation $(2.2), (2.3) $ whenever $X$ does.
 \end{proof}

 \begin{theorem}
 ${\mathcal N} $  has a natural symplectic structure and an almost
 complex structure  compatible with the symplectic form $\Omega $
 and the metric $g$. \label{alcom}
 \end{theorem}

 \begin{proof}
 First we show that the almost complex structure descends to
 ${\mathcal N}$. Then using this and the symplectic quotient
 construction we will show that $\Omega$ gives a symplectic
 structure on ${\mathcal N}$.

 (a) To show that ${\mathcal I}$ descends as  an almost complex
 structure we let $pr: {\mathcal S} \rightarrow {\mathcal S}/G =
 {\mathcal N}$ be the projection map and set $[p] = pr (p)$. Then
 we can naturally identify $T_{[p]}  {\mathcal N} $ with the
 quotient space $T_p {\mathcal S} / T_p O_p, $ where $ O_p = G
 \cdot p $ is the gauge orbit. Using the metric $g$ on ${\mathcal
 S}$ we can realize $T_{[p]} {\mathcal N}$ as a subspace in $T_p
 {\mathcal S}$ orthogonal to  $T_p O_p$. Then by lemma
 ~\ref{ortho}, this subspace is invariant under ${\mathcal I}$.
Thus $I_{[p]} ={\mathcal I} |_{T_p (O_p )^{\perp}}$, gives the desired
almost  complex structure. This construction does not depend on the choice of
$p$ since ${\mathcal I}$ is $G$-invariant.

 (b) The symplectic structure $\Omega$ descends to $\mu^{-1}(0) /
 G$, (by proposition \ref{props} and by the Marsden-Wienstein
 symplectic quotient construction ,~\cite{GS}, ~\cite{H}, since the
 leaves of the characteristic foliation are the gauge orbits). Now,
 as a $2$-form $\Omega$ descends to ${\mathcal N}$,   due to
 proposition (~\ref{inv}) so does the metric $g$. We check that
 equation $(2.2), (2.3),$ does not give rise to new degeneracy of
 $\Omega$ (i.e. the only degeneracy of $\Omega$ is due to $(2.1)$
 but along gauge orbits). Thus $\Omega $ is symplectic on ${\mathcal N}$.
 Since $g$ and ${\mathcal I}$ descend to ${\mathcal N}$ the latter is
 symplectic and almost complex.
 \end{proof}

\subsection{Hyperk\"{a}hler structure in the moduli spaces
${\mathcal M}$ and $\Sigma_{\Psi}$}

We recall that we realised equation $(2.1)$ as a moment map.
To realize the equation $(2.2)$ as a moment map we first rewrite the second
equation as

$(2.2)$ $2\bar{\partial} \Phi^{\prime}  =  -<\psi_1, \psi_2>_H \omega$

where
$\Phi^{\prime} = -i \Phi = -i \phi dz + i \bar{\phi} d \bar{z} \in \Omega^{1}(M, {\RR})$. We rename $\Phi^{\prime}$ as $\Phi$. This notation will be valid only in this section where we do not consider equation $(2.3)$.
We need to define another symplectic form ${\mathcal Q}$ on ${\mathcal C}$,
which is complex-valued,

$ {\mathcal Q} ( X, Y) =   - 2\int_M  \alpha_1 ^{0,1} \wedge \gamma_2 ^{1,0}  + 2\int_M \alpha_2 ^{0,1} \wedge \gamma_1 ^{1,0}
      -  \int_M ( \beta_1 ^1 \bar{ \beta^2 _2 } - \bar{\beta_1 ^2 } \beta_2 ^1 )  \omega $
where $ X  = (\alpha_1, \beta_1, \gamma_1)$, $Y = (\alpha_2, \beta_2, \gamma_2)$ $ \in T_p {\mathcal E}$.

\begin{proposition}
The vector field $X_{\zeta}$ induced by the gauge action is Hamiltonian with
respect to the symplectic form ${\mathcal Q}$.
\end{proposition}

\begin{proof}

Define the Hamiltonian to be $H_{\zeta} (A, \psi, \phi) = \int_M \zeta ( 2\bar{\partial} \Phi +  \psi_1 \bar{\psi_2} \omega ), $ where $\zeta \in \Omega(M, i {\RR })$. Then for $X = (\alpha, \beta_2, \gamma) \in T_p {\mathcal E},$
\begin{eqnarray*}
d H_{\zeta}(X) &=&  \int_M \zeta (2\bar{\partial} \gamma +  ( \beta_2 ^1 \bar{\psi_2} + \psi_1 \bar{\beta_2 ^2} ) \omega ) \\
  &=& 2\int_M -\bar{\partial} \zeta \wedge \gamma  + \int_M   (\zeta \beta^1_2  \bar{\psi_2} +\zeta \psi_1 \bar{\beta}_2 ^2 ) \omega  \\
 &=& 2\int_M -\bar{\partial} \zeta \wedge \gamma  + \int_M  (- \beta_2 ^1 \bar{\zeta} \bar{\psi_2} + \zeta \psi_1 \bar{\beta_2 ^2} ) \omega  \\
 &=& {\mathcal Q}(X_{\zeta}, X).
\end{eqnarray*}
where $X_{\zeta} = (d \zeta, - \zeta \Psi, 0)$.
Thus we can define the moment map of the action with respect to the form
${\mathcal Q}$ to be : $\mu_{{\mathcal Q}} = 2\bar{\partial} \Phi +  <\psi_1 ,\psi_2> \omega.$ Thus  equation (2.2) is precisely $\mu_{{\mathcal Q}} = 0$.
\end{proof}

\begin{proposition}

The configuration space ${\mathcal C}$ has a Riemannian metric

$ g( X, Y) =   \int_M *\alpha_1 \wedge \alpha_2 + \int_M Re< \beta_1, \beta_2> \omega +  \int_M *\gamma_1 \wedge \gamma_2,$

where $X = (\alpha_1,\beta_1 = \left[ \begin{array}{c}
\beta_1 ^1  \\
\beta_1^2
\end{array} \right] ,\gamma_1)$, $Y=(\alpha_2, \beta_2=\left[ \begin{array}{c}
\beta_2 ^1  \\
\beta_2 ^2
\end{array} \right] ,\gamma_2)$
$\in T_p {\mathcal C} $,

and three complex structures

${\mathcal I} = \left[ \begin{array}{ccc}
*  &  0  &   0   \\
0  &  I  & 0    \\
0  &  0  & -*
\end{array} \right] $ , ${\mathcal J}=  \left[ \begin{array}{ccc}
0  &  0  & *  \\
0  &  J  & 0   \\
*  &  0  & 0
\end{array} \right]  $,
${\mathcal K} =  = \left[ \begin{array}{ccc}
0  &  0  & -1  \\
0  &  K  &  0  \\
1  &  0  & 0
\end{array} \right] $,

where $ I = \left[ \begin{array}{cc}
i  & 0 \\
0  & -i
\end{array}
\right]$,   $ J = \left[ \begin{array}{cc}
0 & 1 \\
-1 & 0
\end{array}\right]$ and $ K= \left[ \begin{array}{cc}
0 & i \\
i & 0
\end{array}
\right]$, and $*: \Omega^1(M) \rightarrow \Omega^1(M)$ is the Hodge-star
operator. They satisfy ${\mathcal I} {\mathcal J} = {\mathcal K}$ ,
${\mathcal J} {\mathcal K} = {\mathcal I}$ , and $ {\mathcal K} {\mathcal I} = {\mathcal J}$. The three symplectic structures
$\omega_1(X, Y) = g( {\mathcal I} X, Y )$, $\omega_2 (X,Y) = g ({\mathcal J}X, Y ) $, $\omega_3 (X,Y) = g({\mathcal K} X, Y )$
are such that $\omega_2 + i\omega_3 = {\mathcal Q}$.
\label{hyp}
\end{proposition}

\begin{proof}
\begin{eqnarray*}
 \omega_1(X,Y) &=&  -\int_M \alpha_1 \wedge \alpha_2  +  \int_M Re <I \beta_1, \beta_2 > \omega  +  \int_M \gamma_1 \wedge \gamma_2 \\
\omega_2 (X, Y)  &=&  \int_M -\gamma_1 \wedge \alpha_2  - \int_M \alpha_1 \wedge \gamma_2 + \int_M  Re( \beta_1 ^2 \bar{\beta_2 ^1} - \beta_1 ^1 \bar{ \beta_2 ^2 }) \omega\\
\omega_3 (X, Y) &=& \int_M -*\gamma_1 \wedge \alpha_2 +  \int_M Re< K \beta_1 , \beta_2 > + \int_M  *\alpha_1 \wedge \gamma_2 \\
&=& \int_M -*\gamma_1 \wedge \alpha_2 +   \int_M Re( i\beta_1 ^2 \bar{ \beta_2 ^1 } +i\beta_1 ^1 \bar{ \beta}_2 ^2  ) \omega +\int_M *\alpha_1 \wedge \gamma_2
\end{eqnarray*}
so that indeed
\begin{eqnarray*}
(\omega_2 + i \omega_3 )(X, Y) &=& \int_M (-\gamma_1 - i * \gamma_1 ) \wedge \alpha_2
+   \int_M [Re ( \beta_1 ^2 \bar{ \beta_2 ^1 } - \beta_1 ^1 \bar{\beta_2 ^2 }) \\
&+& i Re(i \beta_1 ^2 \bar{\beta_2 ^1} + i \beta_1 ^1 \bar{ \beta_2 ^2 })]  \omega  + \int (-\alpha_1 + i *\alpha_1) \wedge \gamma_2 \\
&=&  - 2\int_M  (\gamma_1)^{1,0} \wedge \alpha_2 ^{1,0}  -  \int_M ( \beta_1 ^1 \bar{\beta_2 ^2 } - \bar{ \beta_1 ^2} \beta_2 ^1 ) \omega \\
& &  - 2 \int_M (\alpha_1)^{0,1} \wedge \gamma_2 ^{1,0} \\
&=&  {\mathcal Q} (X, Y).\\
\end{eqnarray*}
\end{proof}

Let $\tilde{{\mathcal S}} = \mu^{-1}(0)$ $\cap$   $\mu_{{\mathcal Q}}^{-1}(0) $ $\subset {\mathcal E}$ be the solution space to the equations $(2.1) $ and $(2.2),$ and
denote by ${\mathcal M} = \tilde{{\mathcal S}}/G $  the corresponding moduli space.

\begin{theorem}
Let $M$ be a compact Riemann surface of $g \geq 1$.
Let ${\mathcal M}$ be the moduli space of solutions to equations $(2.1)$ and
$(2.2)$. Then the Riemannian metric $g$ induced by the metric on ${\mathcal C}$ is hyperk\"{a}hlerian, and ${\mathcal M}$ is  hyperk\"{a}hler.
\label{hyper}
\end{theorem}

\begin{proof}
Since ${\mathcal I}, {\mathcal J}, {\mathcal K}$, $g$ and $\omega_1$ ,
${\mathcal Q}$  are $G$-invariant, and ${\mathcal M}$ comes from a symplectic
reduction, it follows that the symplectic forms $\omega_i$, $i = 1,2,3,$
descend to ${\mathcal M}$ as symplectic forms . Also, from  the proof of
theorem (~\ref{alcom}) and proposition (~\ref{hyp}) it follows that
${\mathcal I}, {\mathcal J}, {\mathcal K}$ are well defined almost complex
structures on ${\mathcal M}$. To show that they are integrable, we use the
following lemma of Hitchin (see ~\cite{H}).
\end{proof}
\begin{lemma}
Let $g$ be an almost hyperk\"{a}hler metric, with $2$-forms $\omega_1$,
$\omega_2$, $\omega_3 $ corresponding  to almost complex structures
${\mathcal I}$,$ {\mathcal J}$ and ${\mathcal K}$.
Then $g$ is hyperk\"{a}hler if each $\omega_i$ is closed.
\end{lemma}

\begin{theorem}
Let $M$ be a Riemann surface of genus $g \geq 1$. Fix the equivalence
class of $\Psi$ such that $\psi_1$ and $\psi_2$ are each not identically zero
 and such that $< \psi_1, \psi_2> \omega $ is $\bar{\partial}$-exact.
Then, $\Sigma_{[\Psi]}$ is  hyperK\"{a}hler affine manifold of  dimension
$4g$. Also, ${\mathcal M} = Sp \times   \Sigma_{[\Psi]}$ where
$ Sp = \{ \psi: < \psi_1, \psi_2> \omega $  is  $\bar{\partial}$-exact.$\}$
\label{hyper2}
\end{theorem}

\begin{proof}
On $ {\mathcal W} $ one defines the same symplectic forms $\omega_1$, $\omega_2$and $\omega_3$ as in the previous section. On $\Sigma_{[\Psi]} $ these forms
restrict to

$ \omega_{1 [\Psi]} =  -\int_{M} \alpha^1 \wedge \alpha^2 +   \int_M \gamma ^1 \wedge \gamma^2,$

$ \omega_{2 [\Psi]} = -\int_M \gamma_1 \wedge \alpha_2  - \int_M \alpha_1 \wedge \gamma_2  $

$ \omega_{3 [\Psi]} = -\int_M *\gamma_1 \wedge \alpha_2 +\int_M *\alpha_1 \wedge \gamma_2 $

which are, by arguments same as in the previous section,
 hyperk\"{a}hlerian with respect to the complex structures
${\mathcal I}_1 =  \left[ \begin{array}{cc}
*   & 0   \\
0  & -*
\end{array} \right], $
${\mathcal J}_1 = \left[ \begin{array}{cc}
0   & *  \\
*   & 0                                                                       \end{array} \right]  $, ${\mathcal K}_1 = \left[ \begin{array}{cc}
0   & -1  \\
1   & 0
\end{array} \right] $ and to the Riemannian metric
$g (X, Y) =   \int_M *\alpha_1 \wedge \alpha_2 +  \int_M *\gamma_1 \wedge \gamma_2$
where $X = (\alpha_1, \gamma_1)$ and $Y = (\alpha_2, \gamma_2) $ $\in T_p \Sigma_{[\Psi]}$.
\end{proof}

{\bf Acknowledgements:}

I wish to thank  Professor Leon Takhtajan for  very useful discussions.
I am grateful to Dr. Rajesh Gopakumar, Joe Coffey, Professor Dusa Mc Duff,
Professor Karen Uhlenbeck and Professor Dan Freed   for their interest and
for useful discussions. This work was part of my PhD thesis while I was at
SUNY at Stony Brook, USA and completed in HRI, Allahabad, India.

Harish Chandra Research Institute, Chhatnag, Jhusi, Allahabad, 211019, India.
email: rkmn@mri.ernet.in

\end{document}